**Article**

# Logic and Probability

**Gunn Quznetsov***

## Abstract

The propositional logic is generalized on the real numbers field. The logical function with all properties of the classical probability function is obtained. The logical analog of the Bernoulli independent tests scheme is constructed. The logical analog of the Large Number Law is deduced from properties of these functions.

## *Introduction*

There is the evident nigh affinity between the classical probability function and the Boolean function of the classical propositional logic [1]. These functions are differed by the range of value, only. That is if the range of values of the Boolean function shall be expanded from the two-elements set {0; 1} to the segment [0; 1] of the real numeric axis then the logical analog of the Bernoulli Large Number Law [2] can be deduced from the logical axioms. These topics are considered in this article.

## *Propositional Logic*

**Def. 1:** Sentence «Θ» is *a true sentence* if and only if Θ [3].

For example: sentence «it rains» is the true sentence if and only if it rains.

**Def. 2:** Sentence «Θ» is *a false sentence* if and only if it is not that Θ.

There exist many sentences which are neither true nor false.

**Example 1:** obviously, the sentence "July 23, 2031 in Chelyabinsk is raining" is not true and is not false.

**Example 2:** «the sentence which is written after word «**Example 2:**» is false» - you can test that this sentence is not true and is not false[†].

The sentence from the first example will be true or false after some time. The sentence of the second example, in principle, is neither true nor false. I call sentences which cannot be true and cannot be false *senseless sentences*. And sentences which can be true or false, are *sensible*

* Correspondence: Gunn Quznetsov, Chelyabinsk State University, Chelyabinsk, Russia. E-mail: gunn@mail.ru

† Eubulides of Miletus, Εὐβουλίδης; Liar paradox, fl. 4th century BCE

*sentences*. I denote set of sensible sentences by the symbol A. Therefore, the sentence of first example is element of A, but the sentence of second example is not element of A.

Further we consider only sensible sentences.

**Def. 3:** Sentences *A* and *B* are *equal* (*A* = *B*) if *A* is true if and only if *B* is true.

Hereinafter we use usual notions and assertions of the classical propositional logic (for instance: [4]).

**Def. 4:** The sentence *C* is called *a conjunction* of sentences *A* and *B* (denote: *C* = (*A*&*B*)) if *C* is true, if and only if *A* is true and *B* is true. *A* and *B* are called *conjuncts* of this conjunction.

**Def. 5:** The sentence *C* is called *a negation* of sentences *A* (denote: *C* = (¬*A*)) if *C* is true, if and only if *A* is not true.

## *Classical Propositional Logic*

Let $A_0$ be a set of sentences each of which is either false or true. In this part only elements of $A_0$ are considered.

### Natural Propositional Logic

Further I set out the version of the Gentzen Natural Propositional calculus (NPC) [5]:

Expression «Sentence *C* is a logical consequence of the list of sentences Γ» will be written as the following: «Γ → *C*». Such expressions are called *sequences*. Elements of list Γ are called *hypothesizes*.

**Def. 6**:

1) A sequence of the form *C* → *C* is called NPC-axiom.
2) A sequences of form Γ → *A* and Γ → *B* is obtained from a sequence of form Γ → (*A*&*B*) by *the conjunction removing rule* (denote: R&).
3) A sequence of form $\Gamma_1,\Gamma_2$ → (*A*&*B*) is obtained from a sequence of form $\Gamma_1$ → *A* and a sequence of form $\Gamma_2$ → *B* by *the conjunction inputting rule* (denote: I&).
4) A sequence of form Γ → *C* is obtained from a sequence of form Γ → (¬ (¬*C*)) by *the negation removing* rule (denote: R¬).
5) A sequence of form $\Gamma_1,\Gamma_2$ → (¬*C*) is obtained from a sequence of form $\Gamma_1$,*C* → *A* and from a sequence of form $\Gamma_2$,*C* → (¬*A*) by *the negation inputting rule* (denote: I¬).
6) A finite string of sequences is called *a propositional natural deduction* if every element of this string either is NPC axiom or is received from preceding sequences by one of the deduction rules (R&, I&, R¬, I¬).

These logical rules look naturally in light of the previous definitions. Hence, if a sequence Γ → *A* is contained in some natural propositional deduction, then sentence *A* follows logically from the list of hypotheses Γ.

**Example 3:**

1. $A \to A$, NPC-axiom;
2. $((\neg A)\&(\neg B)) \to ((\neg A)\&(\neg B))$, NPC-axiom;
3. $((\neg A)\&(\neg B)) \to (\neg A)$, R&, 2;
4. $A \to (\neg ((\neg A)\&(\neg B)))$, I¬, 1,3.

This string of sequences is a propositional natural deduction in accordance with point 6 of Def. 6 because every element of this string either is NPC axiom or is received from preceding sequences by one of the deduction rules (R&, I&, R¬, I¬). Since sequence $A \to (\neg ((\neg A)\&(\neg B)))$ is contained in this deduction then sentence $(\neg ((\neg A)\&(\neg B)))$ follows logically from sentence *A*.

**Example 4:**

1. $(A\&(\neg A)) \to (A\&(\neg A))$, NPC-axiom;
2. $(A\&(\neg A)) \to A$, R&, 1;
3. $(A\&(\neg A)) \to (\neg A)$, R&, 1;
4. $\to (\neg (A\&(\neg A)))$, I¬, 2,3.

This string is a propositional natural deduction, too. There sentence $(\neg (A\&(\neg A)))$ follows logically from the empty list of hypothesizes. Such sentences are called *propositionally provable sentences.*

**Boolean functions**

**Def. 7:** Let function g has the double-elements set {0; 1} as a range of reference and $A_0$ as a domain. And let

1) $g(\neg A) = 1 - g(A)$ for every sentence *A*;
2) $g(A\&B) = g(A) \times g(B)$ for all sentences *A* and *B*;

In this case function g is *a Boolean function*.

Hence if g is a Boolean function then for every sentence *A*: $(g(A))^2 = g(A)$ .

A Boolean function can be defined by the following table:

| *A* | *B* | (¬*A*) | (*A*&*B*) |
|---|---|---|---|
| 0 | 0 | 1 | 0 |
| 0 | 1 | 1 | 0 |
| 1 | 0 | 0 | 0 |
| 1 | 1 | 0 | 1 |

Such tables can be constructed for any sentence. For example:

| *A* | *B* | *C* | (¬ ((¬(*A*&(¬*C*))) &((*A*&*B*) &(¬ *C*)))) |
|---|---|---|---|
| 0 | 0 | 0 | 1 |
| 0 | 0 | 1 | 1 |
| 0 | 1 | 0 | 1 |
| 0 | 1 | 1 | 1 |
| 1 | 0 | 0 | 1 |
| 1 | 0 | 1 | 1 |
| 1 | 1 | 0 | 1 |
| 1 | 1 | 1 | 1 |

If $g$ is a Boolean function then by Def.7:

$g(\neg((\neg(A\&(\neg C)))\&((A\&B)\&(\neg C)))) = 1 - g((\neg(A\&(\neg C)))\&((A\&B)\&(\neg C))) =$
$= 1 - g(\neg(A\&(\neg C)))\times g((A\&B)\&(\neg C)) = 1 - (1 - g(A\&(\neg C)))\times g(A\&B)\times g(\neg C) =$
$= 1 - (1 - (g(A)\times g(\neg C)))\times g(A\&B)\times g(\neg C) =$
$= 1 - (g(A\&B)\times g(\neg C) - g(A)\times g(\neg C)\times g(A\&B)\times g(\neg C)) =$
$= 1 - (g(A\&B)\times g(\neg C) - g(A)\times g(A\&B)\times g(\neg C)) =$
$= 1 - (g(A)\times g(B)\times g(\neg C) - g(A)\times g(A)\times g(B)\times g(\neg C)) =$
$= 1 - (g(A)\times g(B)\times g(\neg C) - g(A)\times g(B)\times g(\neg C)) = 1.$

Therefore, for every Boolean function $g$:

$g(\neg((\neg(A\&(\neg C)))\&((A\&B)\&(\neg C)))) = 1.$

Such sentences are called *tautologies*.

**Def. 8:** A set $\mathsf{A}_{0,0}$ of sentences is called *a basic set* if for every element $A$ of $\mathsf{A}_{0,0}$ there exist Boolean functions $g_1$ and $g_2$ such that the following conditions are fulfill:

1) $g_1(A) \neq g_2(A)$;

2) for every element $B$ of set $\mathsf{A}_{0,0}$: if $B \neq A$ then $g_1(B) = g_2(B)$.

Set $\mathsf{A}_{0,0}$ does not contain conjunctions and negations of this set elements because if $(A\&B) \in \mathsf{A}_{0,0}$, $A \in \mathsf{A}_{0,0}$, and $B \in \mathsf{A}_{0,0}$ then Boolean functions $g_1$ and $g_2$ exist such that

$g_1(A\&B) = 0$, $g_2(A\&B) = 1$,
$g_1(A) = g_2(A)$ ,
$g_1(B) = g_2(B)$ .

But it is impossible. Similar argumentation is for and negations.

**Def. 9**: A set $[\mathsf{A}_{0,0}]$ of sentences is called *a propositional closure* of the set $\mathsf{A}_{0,0}$ if the following conditions are satisfied:

1) if $A \in \mathsf{A}_{0,0}$ then $A \in [\mathsf{A}_{0,0}]$;

2) if $A \in [\mathsf{A}_{0,0}]$ then $(\neg A) \in [\mathsf{A}_{0,0}]$;

3) if $A \in [\mathsf{A}_{0,0}]$ and $B \in [\mathsf{A}_{0,0}]$ then $(A\&B) \in [\mathsf{A}_{0,0}]$;

4) there are no other elements of the set $[\mathsf{A}_{0,0}]$ except the enumerated above.

Henceforth, $[\mathsf{A}_{0,0}] = \mathsf{A}_0$.

**Th. 1:** Each naturally propositionally proven sentence is a tautology[‡].

**Th. 2 (Laszlo Kalmar):**[4] Each tautology is a naturally propositionally proven sentence.

***Consequently, whole propositional logic is defined by a Boolean function.***

[‡] Please see the proofs in Appendix

**Th. 3:** Each naturally propositionally proven sentence is a true sentence.

**Th. 4:** Each tautology is the true sentence.

## *Probability*

Further we consider set A (the set of all sensible sentences).

### Events

**Def. 10**: A set B of sentences is called *event, expressed by sentence C*, if the following conditions are fulfilled:

1. $C \in$ B;
2. if $A \in$ B and $D \in$ B then $A = D$;
3. if $D \in$ B and $A = D$ then $A \in$ B.

In this case denote: B := °*C*.

**Def. 11**: An event B *occurs* if here exists a true sentence *A* such that $A \in$ B.

**Def. 12**: Events A and B *equal* (denote: A = B) if A occurs if and only if B occurs.

**Def. 13**: Event C is called *product* of event A and event B (denote: C = (A· B)) if C occurs if and only if A occurs and B occurs.

**Def. 14**: Events C is called *complement* of event A (denote: C = (#A)) if C occurs if and only if A does not occur.

**Def. 15**: (A+ B):= (#((#A)·(#B))). Event (A+ B) is called *sum* of event A and event B.

Therefore, the sum of event occurs if and only if there is at least one of the addends.

**Def. 16**: *The persistent event* (denote: T) is the event which contains a tautology.

Hence, T occurs by Th.4.

**Def. 17**: *The impossible event* (denote: F) is event which contains negation of a tautology.

Hence, F does not occur by Th.4, too.

### B-functions

**Def. 18**: Let ɓ(X) be any function defined on the set of events.

And let the real numbers segment [0; 1] is this function frame reference.
Let there exists an event $C_0$ such that ɓ($C_0$) = 1.
Let for all events *A* and *B*:

ɓ(A·B) + ɓ(A·(#B)) = ɓ(A).

In that case function ɓ(X) is called *a B-function*.
By this definition:

ɓ(A·B) ≤ ɓ(A). (p1)

Hence, ɓ(T·$C_0$) ≤ ɓ(T). Because T·$C_0$ = $C_0$ (by Def.13 and Def.16) then ɓ($C_0$) ≤ ɓ(T). Because ɓ($C_0$) = 1 then

ɓ(T) = 1. (p2)

From Def.18: ɓ(T·B) + ɓ(T·(#B)) = ɓ(T). Because T·D = D for any D by Def.13 and Def.16 then ɓ(B) + ɓ (#B) = ɓ(T). Hence, by (p2): for any B:

ɓ(B) + ɓ (#B) =1. (p3)

Therefore, ɓ(T) + ɓ (#T) =1. Hence, in accordance (p2) and in accordance Def.14, Def.16, and Def.17 : 1 + ɓ (F) =1. Therefore,

ɓ (F) = 0. (p4)

In accordance with Def.18, Def.15, and (p3):
ɓ(A·(B + C)) = ɓ(A·(#((#B)·(#C)))) = ɓ(A) − ɓ(A·((#B)·(#C))) = ɓ(A) − ɓ((#C)·((#B)· A)) =
= ɓ(A) − (((#B)· A) − ɓ(C·((#B)· A))) = ɓ(A) − ((#B)· A) + ɓ(C·((#B)· A)) =
= ɓ(A·B) + (A·C) − ɓ(A·B·C)).

And

ɓ( (A·B) + (A·C)) = ɓ(#((#(A·B))·(#(A·C)))) = 1 − ɓ((#(A·B))·(#(A·C))) =
= 1 − (1 − ɓ(A·B)) + (ɓ(A·C) − ɓ(A·B·A·C)) = ɓ(A·B)) + ɓ(A·C) − ɓ(A·B·A·C) =
= ɓ(A·B)) + ɓ(A·C) − ɓ(A·B·C)

because A·A = A in accordance with Def.13.

Therefore:

ɓ(A·(B +C)) = ɓ(A·B) + (A·C) − ɓ(A·B·C)) (p5)

and

ɓ( (A·B) + (A·C)) = ɓ(A·B)) + ɓ(A·C) − ɓ(A·B·C). (p6)

Hence (distributivity):

ɓ(A·(B + C)) = ɓ( (A·B) + (A·C)). (p7)

If A = T then from (p5) and (p6) (*the addition formula of probabilities*):

ɓ(B+C) = ɓ(B) + ɓ(C) − ɓ(B·C). (p8)

**Def. 19**: Events B and C are *antithetical* events if (B· C) = F.
From (p8) and (p4) for antithetical events B and C:

ɓ(B + C) = ɓ(B) + ɓ(C). (p9)

**Def. 20**: Events B and C are *independent for B-function* ɓ events if ɓ(B·C)= ɓ(B)·ɓ(B).
If events B and C are independent for B-function ɓ events then:

ɓ(B·(#C)) = ɓ(B) − ɓ(B·C) = ɓ(B) − ɓ(B)·ɓ(C) = ɓ(B)·(1 − ɓ(C)) = ɓ(B)·ɓ(#C).

Hence, if events B and C are independent for B-function b events then:

ɓ(B·(#C)) = ɓ(B)·ɓ(#C). (p10)

Let calculate:

ɓ(A·(#A)·C)=ɓ(A·C) − ɓ(A·A·C) = ɓ(A·C) − ɓ(A·C) = 0. (p11)

## Independent Tests

Let **N** be the natural numbers set.

**Def. 21:** Let st(*n*) be a function such that st(*n*) has domain on **N** and has a range of values in the set of events.

In this case an event C is *a* [st]-*series of range r with V-number k* if C, *r* and *k* is subject to one of the following conditions:

1) *r* =1 and *k* = 1, C := st(1), or *k* = 0, C := (#st(1));
2) B is a [st]-series of range *r* − 1 with V-number *k* − 1 and

C:= (B· st(*r*)) ,

or B is a [st]-series of range *r* − 1 with V-number *k* and

C:= (B·(#st(*r*))).

Let us denote a set of [st]-series of range *r* with V-number *k* as [st](*r; k*).

For example, if st(*n*) is event $B_n$ then the following events:

($B_1$·$B_2$·(#$B_3$)), ($B_1$·(#$B_2$)·$B_3$), ((#$B_1$)·$B_2$·$B_3$)

are elements of [st](3;2), and

$(B_1 \cdot B_2 \cdot (\#B_3) \cdot B_4 \cdot (\#B_5)) \in [st](5;3)$.

**Def. 22:** Def. 4.2.2: Function st($n$) is *independent for B-function* ե if:

ե(st(1)· st(2)·… ·st($k$)) = ե(st(1))·ե(st(2))·… ·ե(st($k$))
for any $k$.

**Def. 23:** Let st($n$) has domain on the set of natural numbers and has range of values in the set of events.

In this case event C is called *a* [st]-*sum of range r with V-number k* (denote: C:= ŧ[st]($r$, $k$)) if C is a sum of all elements of [st]($r$, $k$).

For example, if st($n$) is the sentence $C_n$ then:

$((\#C_1) \cdot (\#C_2) \cdot (\#C_3))$ = ŧ[st] (3;0),

ŧ[st] (3;2) = $(((\#C_1) \cdot C_2 \cdot C_3) + (C_1 \cdot (\#C_2) \cdot C_3) + (C_1 \cdot C_2 \cdot (\#C_3))$,

ŧ[st] (3;1) = $((C_1 \cdot (\#C_2) \cdot (\#C_3)) + ((\#C_1) \cdot C_2 \cdot (\#C_3)) + ((\#C_1) \cdot (\#C_2) \cdot C_3))$,

$(C_1 \cdot C_2 \cdot C_3)$ = ŧ[st] (3;3).

**Def. 24**: Let a function $s_A(n)$ be defined on **N**, has range of values in the set of events, and be independent for a B-function ե.

And let $s_A(n)$ satisfies the following condition: ե($s_A(n)$) = ե(A) for any $n$.

In that case the [$s_A$]-series of rank $r$ with V-number $k$ is called *series of r independent for B-function* ե [$s_A$]-*tests of event* A *with result k.*

**Def. 25:** Function $\nu_r[s_A]$ is called *a frequency of event A in* [$s_A$]-*series* if $\nu_r[s_A] = k/r$ if and only if event ŧ[$s_A$]($r$, $k$) occurs.

Hence,

°«$\nu_r(s_A) = k/r$» = ŧ[$s_A$]($r$, $k$). (p12)

**Th. 5:** (**the Bernoulli Formula**) [2] If s($n$) is independent for B-function ե and there exists a real number $p$ such that for all $n$: ե(s($n$)) = $p$ then

$$\mathrm{b}(ŧ[s](r,k)) = \frac{r!}{k!\,(r-k)!} p^k (1-p)^{r-k}.$$

**Def. 26:** Let a function s($n$) be defined on **N** and has a range of values in the set of events.

In that case an event Ŧ[s]($r$,$k$,$l$) with natural $r$, $k$, $l$ is defined in the following way:

1) Ŧ[s]($r$,$k$,$k$) := ŧ[s]($r$,$k$),

2) Ŧ[s]($r$,$k$,$l$ + 1) := (Ŧ[s]($r$,$k$,$l$) + ŧ[s]($r$,$l$ + 1)).

If $a$ and $b$ are real numbers, and $k - 1 < a \leq k$ and $l \leq b < l + 1$ then Ŧ[s]($r$,$a$,$b$) := Ŧ[s]($r$,$k$,$l$).

**Th. 6:** Ŧ[$s_A$]($r,a,b$) occurs if and only if $a/r \leq \nu_r[s_A] \leq b/r$.

**Th.7:** If s($n$) is independent for a B-function ɓ and there exists a real number $p$ such that ɓ(s($n$)) = $p$ for all $n$ then

$$\mathrm{b}\left(Ŧ[s_A](r,a,b)\right) = \sum_{a \leq k \leq b} \frac{r!}{k!\,(r-k)!} p^k (1-p)^{r-k}.$$

**Th. 8:** If s($n$) is independent for a B-function ɓ and there exists a real number $p$ such that ɓ(s ($r$)) = $p$ for all $r$ then

$$\mathrm{b}\left(Ŧ[s_A]\left(r, r\cdot(p-\varepsilon), r\cdot(p+\varepsilon)\right)\right) \geq 1 - \frac{p\cdot(1-p)}{r\cdot\varepsilon^2}$$

for every positive real number $\varepsilon$ .

Hence, in accordance with Th.6:

$$\mathrm{b}(\text{°"}p-\varepsilon \leq \nu_r[s_A] \leq p+\varepsilon\text{"}) \geq 1 - \frac{p\cdot(1-p)}{r\cdot\varepsilon^2}.$$

The right part of this inequality doesn't depend on sequence s. Hence it can be rewritted as the following:

$$\mathrm{b}(\text{°"}p-\varepsilon \leq \nu_r[\mathrm{A}] \leq p+\varepsilon\text{"}) \geq 1 - \frac{p\cdot(1-p)}{r\cdot\varepsilon^2}. \qquad \text{(p13)}$$

## Function of Probability

Nonstandard Numbers

Further some variant of the Robinson non-standard analysis (for instant [6]) is required:

**Def. 27:** A *n-part-set* **S** of **N** is defined recursively as follows:

1) $\mathbf{S}_1 = \{1\}$;
2) $S_{(n+1)} = \mathbf{S}_n \cup \{n+1\}$.

**Def. 28:** If $\mathbf{S}_n$ is a $n$-part-set of **N** and $\mathbf{A} \subseteq \mathbf{N}$ then $\|\mathbf{A}\cap\mathbf{S}_n\|$ is quantity of elements of set $\mathbf{A}\cap\mathbf{S}_n$, and if $\varpi_n(\mathbf{A}) := \|\mathbf{A}\cap\mathbf{S}_n\|/n$ then $\varpi_n(\mathbf{A})$ is called *a frequency* of set **A** on the $n$-part-set $\mathbf{S}_n$.

Because $\varpi_n(\mathbf{N}) = \|\mathbf{N}\cap\mathbf{S}_n\|/n = n/n$ then

$$\varpi_n(\mathbf{N}) = 1. \qquad \text{(s1)}$$

Becaus $\varpi_n(\mathbf{A}\cap\mathbf{B}) + \varpi_n((\mathbf{N}-\mathbf{A})\cap\mathbf{B}) = \|(\mathbf{A}\cap\mathbf{B})\cap\mathbf{S}_n\|/n +$
$+ \|((\mathbf{N}-\mathbf{A})\cap\mathbf{B})\cap\mathbf{S}_n\|/n$ then

$\varpi_n(\mathbf{A} \cap \mathbf{B}) + \varpi_n\left((\mathbf{N} - \mathbf{A}) \cap \mathbf{B}\right) = \varpi_n(\mathbf{B}).$ (s2)

Hence,

$\varpi_n(\mathbf{A} \cap \mathbf{N}) + \varpi_n((\mathbf{N} - \mathbf{A}) \cap \mathbf{N}) = \varpi_n(\mathbf{N})$ and because for any $\mathbf{A}$: $(\mathbf{A} \cap \mathbf{N}) = \mathbf{A}$ then

$\varpi_n(\mathbf{A}) + \varpi_n(\mathbf{N} - \mathbf{A}) = 1.$ (s3)

Therefore, $\varpi_n(\mathbf{N}) + \varpi_n(\mathbf{N} - \mathbf{N}) = 1$. That is $\varpi_n(\mathbf{N}) + \varpi_n(\varnothing) = 1$.

Hence,
$\varpi_n(\varnothing) = 0.$ (s4)

**Def.29:** If "lim" is the Cauchy-Weierstrass "limit" then:

$$\Phi ix := \left\{\mathbf{A} \subseteq \mathbf{N} \mid \lim_{n \to \infty} \varpi_n(\mathbf{A}) = 1\right\}.$$

Hence, in accordance with (s1)

$\mathbf{N} \in \Phi ix$ and $\varnothing \notin \Phi ix$. (s5)

If $\mathbf{B} \in \Phi ix$ then $\lim_{n\to\infty} \varpi_n(\mathbf{N} - \mathbf{B}) = 0$. In accordance with (s2):

$\varpi_n(\mathbf{A} \cap (\mathbf{N} - \mathbf{B})) \leq \varpi_n(\mathbf{N} - \mathbf{B})$. Therefore, $\lim_{n\to\infty} \varpi_n(\mathbf{A} \cap (\mathbf{N} - \mathbf{B})) = 0$. Hence,

$\lim_{n\to\infty} \varpi_n(\mathbf{A} \cap \mathbf{B}) = \lim_{n\to\infty} \varpi_n(\mathbf{A})$.

Therefore, if $\mathbf{B} \in \Phi ix$ and $\mathbf{A} \in \Phi ix$ then $(\mathbf{A} \cap \mathbf{B}) \in \Phi ix$. (s6)

Moreover,

if $\mathbf{A} \in \Phi ix$ and $\mathbf{A} \subseteq \mathbf{B}$ then $\mathbf{B} \in \Phi ix$. (s7)

Therefore, in accordance with (s5), (s6), (s7), $\Phi ix$ is *a filter* (for instance, [6], p.45), but $\Phi ix$ is not an ultrafilter because there exist subsets $\mathbf{A}$ of $\mathbf{N}$ such that $\mathbf{A} \notin \Phi ix$ and $(\mathbf{N} - \mathbf{A}) \notin \Phi ix$.

**Def. 30**: A series of real numbers $\langle r_n \rangle$ and $\langle s_n \rangle$ are *Q-equivalent* (denote: $\langle r_n \rangle \sim \langle s_n \rangle$) if

$\{n \in \mathbf{N} \mid r_n = s_n\} \in \Phi ix.$

Hence, if r, s,u are series of real numbers then r ~ r; if r ~ s then s ~ r; and if r ~ s, and s ~u then r ~u. Therefore, «~» is *an equivalence relation*.

**Def. 31:** *A Q-number* is a set of Q-equivalent series of real numbers.

That is if ***a*** is a Q-number and r$\in$ ***a*** and s $\in$ ***a*** then r ~ s; and if r$\in$ ***a***, and r ~ s then s $\in$ ***a***.

**Def. 32:** A Q-number ***b*** is *a standard Q-number b* if *b* is some real number and there exists a series $\langle r_n \rangle$ such that $\langle r_n \rangle \in$ ***b***, and

$$\{n \in \mathbf{N} \mid r_n = b\} \in \Phi ix.$$

In this case ***b*** := *b*.

**Def. 33:** Q-numbers ***a*** and ***b*** *equal* (denote: ***a*** = ***b***) if ***a*** $\subseteq$ ***b*** and ***b*** $\subseteq$ ***a***.

**Def. 34:** Q-number ***c*** is *sum* of Q-number ***a*** and Q-number *b* (denote: ***c*** = ***a*** + ***b***) if there exist series of real numbers $\langle r_n \rangle$, $\langle s_n \rangle$, $\langle u_n \rangle$ such that $\langle r_n \rangle \in$ ***a***, $\langle s_n \rangle \in$ ***b***, $\langle u_n \rangle \in$ ***c***, and

$$\{n \in \mathbf{N} \mid r_n + s_n = u_n\} \in \Phi ix.$$

If *a* is a real number then *a* + ***b*** = ***a*** + ***b*** where ***a*** is standard Q-number *a*.

**Def. 35:** Q-number ***c*** is *product* of Q-number ***a*** and Q-number ***b*** (denote: ***c*** = ***a*** · ***b***) if there exist series of real numbers $\langle r_n \rangle$, $\langle s_n \rangle$, $\langle u_n \rangle$ such that $\langle r_n \rangle \in$ ***a***, $\langle s_n \rangle \in$ ***b***, $\langle u_n \rangle \in$ ***c***, and

$$\{n \in \mathbf{N} \mid r_n \cdot s_n = u_n\} \in \Phi ix.$$

Hence, ***a*** – ***b*** = ***a*** + (–1) · ***b*** = ***a*** + (***-1***) · ***b***. And

$$"\boldsymbol{c} = \frac{\boldsymbol{a}}{\boldsymbol{b}}" := "\boldsymbol{a} = \boldsymbol{b} \cdot \boldsymbol{c}".$$

**Def. 36:** A Q-number ***x*** is called *an infinitesimal Q-number* if there exists a series of real numbers $\langle x_n \rangle$ such that $\langle x_n \rangle \in$ ***x***, and for all natural numbers *m*:

$$\{n \in N \mid |x_n| < \frac{1}{m}\} \in \Phi ix.$$

Denote by **I** the set of all infinitesimal Q-numbers.

**Def. 37:** Q-numbers ***x*** and ***y*** are *infinitely near* (denote: ***x*** ≈ ***y***) if either (***x*** – ***y***) = 0 or (***x*** – ***y***)$\in$ **I**.

**Def. 38:** A Q-number ***x*** is called *an infinitely large Q-number* if there exists a series $\langle r_n \rangle$ of real numbers such that $\langle r_n \rangle \in$ ***x***, and for every natural number *m*:

$\{n \in \mathbf{N} \mid m < r_n\} \in \mathbf{\Phi}ix.$

Let $\boldsymbol{n}$ be the Q-number which contains the following series

$\langle n \rangle := 1,2,3,4,\ldots,n,\ldots$ .

Let $m$ be some natural number.
In that case:

$$\lim_{n\to\infty} \varpi_n(\{n \in \mathbf{N} \mid m < n\}) = \lim_{n\to\infty} \frac{n-m}{n} = 1.$$

Hence, for any natural $m$:
$\{n \in \mathbf{N} \mid m < n\} \in \mathbf{\Phi}ix.$

Therefore, $\boldsymbol{n}$ is an infinitely large Q-number. Denote $\underline{n}$ *the natural infinity*.

Let $a$ be a positive real number. In this case $a/\boldsymbol{n}$ contains the series $\langle a/n \rangle$. Let $m$ be some natural number and let $k$ be some natural number which is more than $a$. In that case if $n > mk$ then $(a/n) < 1/m$. That is for any natural number $m$:

$$\lim_{n\to\infty} \varpi_n\left(\{n \in \mathbf{N} \mid \frac{a}{n} < \frac{1}{m}\}\right) = \lim_{n\to\infty} \frac{n-mk}{n} = 1. \qquad \text{(s8)}$$

Therefore, $a/\boldsymbol{n}$ is an infinitesimal Q-number in accordance with Def.36.

**Def. 39:** Let A($x$) be a sentence which contains a real number $x$. And let $\boldsymbol{r}$ be a Q-number. In that case event °A($\boldsymbol{r}$) occurs if and only if here a series $\langle r_n \rangle$ of real number exists for which the following conditions are fulfilled: $\langle r_n \rangle \in \boldsymbol{r}$ and

$\{n \in \mathbf{N} \mid {}^\circ\mathrm{A}(r_n) \text{ occurs}\} \in \mathbf{\Phi}ix.$

## P-functions

**Def. 40:** A B-function $\mathcal{P}$ is called *P-function* if for every event A the following condition is fulfilled:

If $\mathcal{P}$(A) ≈ 1 then A occurs.

In accordance with (p13): for any natural number $n$ and for positive real $\varepsilon$:

$$\mathcal{P}({}^{\circ}\text{"}p - \varepsilon \leq \nu_n[\mathrm{A}] \leq p + \varepsilon\text{"}) \geq 1 - \frac{p \cdot (1-p)}{n \cdot \varepsilon^2}.$$

Hence,

$$\mathcal{P}({}^{\circ}\text{"}p - \varepsilon \leq \nu_{\boldsymbol{n}}[\mathrm{A}] \leq p + \varepsilon\text{"}) \geq 1 - \frac{p \cdot (1-p)}{\boldsymbol{n} \cdot \varepsilon^2}.$$

Because in accordance with (s8) $((p \cdot (1-p))/(\boldsymbol{n} \cdot \varepsilon^2)) \in \mathbf{I}$ then in accordanse with Def.37:

$$\mathcal{P}({}^{\circ}\text{"}p - \varepsilon \leq \nu_{\boldsymbol{n}}[\mathrm{A}] \leq p + \varepsilon\text{"}) \approx 1.$$

Hence, event "$p - \varepsilon \leq \nu_n[\mathrm{A}] \leq p + \varepsilon$" occurs.

Since $p = \mathcal{P}(\mathrm{A})$ then for all arbitrarily small real positive $\varepsilon$:

$$\mathcal{P}(\mathrm{A}) = \nu_n[\mathrm{A}] \pm \varepsilon.$$

Consequently, this function has a statistical meaning. Therefore, ***in all over the world there exists the only single such function because values of this function can be defined by repetition of independent tests experimentally***. Therefore, I call this function *the probability function* (proof of of the consistency see in [7]).

**Probability and Logic**

Let $\mathcal{P}$ be the probability function and let $\mathfrak{B}$ be the set of events $\mathrm{A}$ such that either $\mathrm{A}$ occurs or (#$\mathrm{A}$) occurs.

In this case if $\mathcal{P}(\mathrm{A}) = 1$ then $\mathrm{A}$ occurs, and $(\mathrm{A} \cdot \mathrm{B}) = \mathrm{B}$ in accordance with Def.13. Consequently, if $\mathcal{P}(\mathrm{B}) = 1$ then $\mathcal{P}(\mathrm{A} \cdot \mathrm{B}) = 1$. Hence, in this case $\mathcal{P}(\mathrm{A} \cdot \mathrm{B}) = \mathcal{P}(\mathrm{A}) \cdot \mathcal{P}(\mathrm{B})$.

If $\mathcal{P}(\mathrm{A}) = 0$ then $\mathcal{P}(\mathrm{A} \cdot \mathrm{B}) = \mathcal{P}(\mathrm{A}) \cdot \mathcal{P}(\mathrm{B})$ because $\mathcal{P}(\mathrm{A} \cdot \mathrm{B}) \leq \mathcal{P}(\mathrm{A})$ in accordance with (p1).

Moreover in accordance with (p3): $\mathcal{P}(\#\mathrm{A}) = 1 - \mathcal{P}(\mathrm{A})$ since the function $\mathcal{P}$ is a B-function.

If event $\mathrm{A}$ occurs then $(\mathrm{A} \cdot \mathrm{B}) = \mathrm{B}$ and $(\mathrm{A} \cdot (\#\mathrm{B})) = (\#\mathrm{B})$. Hence,

$$\mathcal{P}(\mathrm{A} \cdot \mathrm{B}) + \mathcal{P}(\mathrm{A} \cdot (\#\mathrm{B})) = \mathcal{P}(\mathrm{A}) = \mathcal{P}(\mathrm{B}) + \mathcal{P}(\#\mathrm{B}) = 1.$$

Consequently, if an element $\mathrm{A}$ of $\mathfrak{B}$ occurs then $\mathcal{P}(\mathrm{A}) = 1$. If $\mathrm{A}$ does not occurs then (#$\mathrm{A}$) occurs. Hence, $\mathcal{P}(\#\mathrm{A}) = 1$ and because $\mathcal{P}(\mathrm{A}) + \mathcal{P}(\#\mathrm{A}) = 1$ then $\mathcal{P}(\mathrm{A}) = 0$.

Therefore, on $\mathfrak{B}$ the range of values of $\mathcal{P}$ is the two-element set {0;1} similar the Boolean function range of values.

Hence, on set $\mathfrak{B}$ the probability function obeys definition of a Boolean function (Def.7).

## *Conclusion*

Therefore, ***the probability is logic of the events which have not occurred yet.***

## *Appendix*

**Lm. 1:** If $g$ is a Boolean function then every natural propositional deduction of sequence $\Gamma \rightarrow A$ satisfy the following condition: if $g(A) = 0$ then there exists a sentence $C$ such that $C \in \Gamma$ and $g(C) = 0$.

**Proof of Lm. 1**: is realized by a recursion on number of sequences in the deduction of $\Gamma \vdash A$:

**1. Basis of recursion:** Let the deduction of $\Gamma \rightarrow A$ contains 1 sequence.

In that case a form of this sequence is $A \rightarrow A$ in accordance with the propositional natural deduction definition (Def. 6). Hence in this case the lemma holds true.

**2. Step of recursion:** *The recursion assumption:* Let the lemma holds true for every deduction containing no more than *n* sequences.

Let the deduction of $\Gamma \to A$ contains $n + 1$ sequences.

In that case either this sequence is a NPC-axiom or $\Gamma \to A$ is obtained from previous sequences by one of deduction rules.

If $\Gamma \to A$ is a NPC-axiom then the proof is the same as for the recursion basis.

a) Let $\Gamma \to A$ be obtained from a previous sequence by *R*&.

In that case a form of this previous sequence is either the following $\Gamma \to (A\&B)$ or is the following $\Gamma \to (B\&A)$ in accordance with the definition of deduction. The deduction of this sequence contains no more than *n* elements. Hence the lemma holds true for this deduction in accordance with the recursion assumption. If g $(A) = 0$ then g $(A\&B) = 0$ and g $(B\&A) = 0$ in accordance with the Boolean function definition (Def. 2.10). Hence there exists a sentence *C* such that $C \in \Gamma$ and $g(C) = 0$ in accordance with the lemma.

Hence in that case the lemma holds true for the deduction of sequence $\Gamma \to A$.

b) Let $\Gamma \to A$ be obtained from previous sequences by I&.

In that case forms of these previous sequences are $\Gamma_1 \to B$ and $\Gamma_2 \to G$ with $\Gamma = \Gamma_1,\Gamma_2$ and $A = (B\&G)$ in accordance with the definition of deduction. The lemma holds true for deductions of sequences $\Gamma_1 \to B$ and $\Gamma_2 \to G$ in accordance with the recursion assumption because these deductions contain no more than *n* elements. In that case if $g(A) = 0$ then $g(B) = 0$ or $g(G) = 0$ in accordance with the Boolean function definition. Hence there exist a sentence *C* such that $g(C) = 0$ and $C \in \Gamma_1$ or $C \in \Gamma_2$.

Hence in that case the lemma holds true for the deduction of sequence $\Gamma \to A$.

c) Let $\Gamma \to A$ be obtained from a previous sequence by R¬.

In that case a form of this previous sequence is the following: $\Gamma \to (\neg (\neg A))$ in accordance with the definition of deduction. The lemma holds true for the deduction of this sequence in accordance with the recursion assumption because this deduction contains no more than *n* elements. If $g(A) = 0$ then $g(\neg (\neg A)) = 0$ in accordance with the Boolean function definition. Hence there exists a sentence *C* such that $C \in \Gamma$ and $g(C) = 0$.

Hence the lemma holds true for the deduction of sequence $\Gamma \to A$.

d) Let $\Gamma \to A$ be obtained from previous sequences by I¬.

In that case forms of these previous sequences are $\Gamma_1,G \to B$ and $\Gamma_2,G \to (\neg B)$ with $\Gamma = \Gamma_1, \Gamma 2$, and $A = (\neg G)$ in accordance with the definition of deduction. The lemma holds true for the deductions of sequences $\Gamma_1,G \to B$ and $\Gamma_2,G \to (\neg B)$ in accordance with the recursion assumption because these deductions contain no more than *n* elements.

If g $(A) = 0$ then g $(G) = 1$ in accordance with the Boolean function definition.

Either g $(B) = 0$ or g $(\neg B) = 0$ by the same definition. Hence there exists a sentence *C* such that either $C \in \Gamma_1,G$ or $C \in \Gamma_2,G$ and $g(C) = 0$ in accordance with the recursion assumption.

Hence in that case the lemma holds true for the deduction of sequence $\Gamma \to A$.

**The recursion step conclusion:** Therefore, in each possible case, if the lemma holds true for a deduction contained no more than *n* elements then the lemma holds true for a deduction contained $n + 1$ elements.

**The recursion conclusion:** Therefore the lemma holds true for a deduction of any length.

**Proof of Th. 1:** If a sentence $A$ is naturally propositionally proven then there exists a natural propositional deduction of form $\rightarrow A$. Hence, for every Boolean function $g$: $g(A) = 1$ in accordance with Lm.1. Hence, sentence $A$ is a tautology.

**Designation 1:** Let $g$ be a Boolean function. In that case for every sentence $A$:

$$A^g := \begin{cases} A, if\ g(A) = 1; \\ (\neg A), if\ g(A) = 0. \end{cases}$$

**Lm. 2:** [4] Let $B_1;B_2,\ldots,B_k$ be elements of a basic set $\mathsf{A}_{0,0}$ making up a sentence $A$ by the logical connectors (¬, &).

Let $g$ be any Boolean function.

In that case there exist a propositional natural deduction of sequence

$B_1^g;B_2^g,\ldots,B_k^g \rightarrow A^g$.

**Proof of Lm. 2:** is realized by a recursion on the number of the logical connectors in sentence $A$.

**1.Basis of recursion**: Let $A$ does not contain the logical connectors.

In this case the following string of one sequence:

1. $A^g \rightarrow A^g$, NPC-axiom.

gives the proof of the lemma**.**

**2. Step of recursion:** *The recursion assumption:* Let the lemma holds true for every sentence, containing no more than $n$ logical connectors.

Let sentence $A$ contains $n + 1$ connector.

Let us consider all possible cases:

a) Let $A = (\neg G)$.

In that case the lemma holds true for $G$ in accordance with the recursion assumption because $G$ contains no more than $n$ connectors. Hence, there exists a deduction of sequence

$$B_1^g;B_2^g,\ldots,B_k^g \rightarrow G^g, \qquad (1)$$

here $B_1;B_2,\ldots,B_k$ are elements of basic set making up sentence $G$.

Hence, $B_1;B_2,\ldots,B_k$ make up sentence $A$.

If g $(A) = 1$ then $A^g = A = (\neg G)$ in accordance with Designation 1.

In that case $g(G) = 0$ in accordance with the Boolean function definition. Hence, $G^g = (\neg G) = A$ in accordance with Designation 1.

Hence, in that case a form of sequence (1) is the following:

$B_1^g;B_2^g,\ldots,B_k^g \rightarrow A^g$.

Hence, in that case the lemma holds true.

If $g(A) = 0$ then $A^g = (\neg A) = (\neg (\neg G))$ .in accordance with Designation 1. In that case **g**$(G) = 1$ in accordance with the Boolean function definition. Hence, $G^g = G$ in accordance with Designation 1.

Hence, in that case a form of sequence (1) is

$B_1{}^g;B_2{}^g,\ldots,B_k{}^g \to G$.

Let us continue the deduction of this sequence in the following way:

1. $B_1{}^g;B_2{}^g,\ldots,B_k{}^g \to G$
2. $(\neg G) \to (\neg G)$, NPC-axiom.
3. $B_1{}^g;B_2{}^g,\ldots,B_k{}^g \to (\neg (\neg G))$, I¬ from 1. and 2.

It is a deduction of sequence

$B_1{}^g;B_2{}^g,\ldots,B_k{}^g \to A^g$.

Hence, in that case the lemma holds true.

b) Let $A = (G\&R)$.

In that case the lemma holds true both for $G$ and for $R$ in accordance with the recursion assumption because $G$ and $R$ contain no more than $n$ connectors. Hence, there exist deductions of sequences

$$B_1{}^g;B_2{}^g,\ldots,B_k{}^g \to G^g \quad (2)$$

and

$$B_1{}^g;B_2{}^g,\ldots,B_k{}^g \to R^g, \quad (3)$$

here $B_1;B_2,\ldots,B_k$ are elements of basic set making up sentences $G$ and $R$. Hence $B_1;B_2,\ldots,B_k$ make up sentence $A$.

If $g(A) = 1$ then $A^g = A = (G\&R)$ in accordance with Designation 1.

In that case $g(G) = 1$ and $g(R) = 1$ in accordance with the Boolean function definition.

Hence, $G^g = G$ and $R^g = R$ in accordance with Designation 1.

Let us continue deductions of sequences (2) and (3) in the following way:

1. $B_1{}^g;B_2{}^g,\ldots,B_k{}^g \to G$, (2).
2. $B_1{}^g;B_2{}^g,\ldots,B_k{}^g \to R^g$, (3).
3. $B_1{}^g;B_2{}^g,\ldots,B_k{}^g \to (G\&R)$, I& from 1. and 2.

It is deduction of sequence
$B_1{}^g;B_2{}^g,\ldots,B_k{}^g \to A^g$.

Hence, in that case the lemma holds true.

If g $(A) = 0$ then $A^g = (\neg A) = (\neg (G\&R))$ in accordance with Designation 1.

In that case $g(G) = 0$ or $g(R) = 0$ in accordance with the Boolean function definition.

Hence, $G^g = (\neg G)$ or $R^g = (\neg R)$ in accordance with Designation 1.

Let $G^g = (\neg G)$.

In that case let us continue a deduction of sequence (2) in the following way:

1. $B_1{}^g;B_2{}^g,\ldots,B_k{}^g \to (\neg G)$, (2).
2. $(G\&R) \to (G\&R)$, NPC-axiom.
3. $(G\&R) \to G$, R& from 2.
4. $B_1{}^g;B_2{}^g,\ldots,B_k{}^g \to (\neg (G\&R))$, I¬ from 1. and 3.

It is a deduction of sequence

$$B_1{}^g;B_2{}^g,\ldots,B_k{}^g \rightarrow A^g.$$

Hence, in that case the lemma holds true.
The same result is received if $R^g = (\neg R)$.

**The recursion step conclusion:** If the lemma holds true for sentences contained no more than $n$ connectors then the lemma holds true for sentences contained $n + 1$ connectors.

**The recursion conclusion:** The lemma holds true for sentences, containing any number connectors.

**Proof of Th. 2:** Let sentence $A$ be a tautology. That is for every Boolean function $g$: $g(A) = 1$.
Hence there exists a deduction for sequence

$$B_1{}^g;B_2{}^g,\ldots,B_k{}^g \rightarrow A \qquad (4)$$

for every Boolean function $g$ in accordance with Lm. 2.
There exist Boolean functions $g_1$ and $g_2$ such that

$g_1(B_1) = 0$, $g_2(B_1) = 1$,
$g_1(B_s) = g_2(B_s)$ for $s \in \{2,\ldots, k\}$.

in accordance with Def. 8 because all $B_s$ ($s \in \{1; 2,\ldots, k\}$) are elements of the basic set.
Forms of sequences (4) for these Boolean functions are the following[§]:

$$(\neg B_1), B_2{}^{g_1},\ldots, B_k{}^{g_1} \vdash A, \qquad (5)$$

$$B_1, B_2{}^{g_1},\ldots, B_k{}^{g_1} \vdash A. \qquad (6)$$

Let us continue deductions of these sequences in the following way:

1. $(\neg B_1), B_2{}^{g_1},\ldots, B_k{}^{g_1} \vdash A$,, (5),
2. $B_1, B_2{}^{g_1},\ldots, B_k{}^{g_1} \vdash A$, (6),
3. $(\neg A) \vdash (\neg A)$, NPC-axiom.
4. $(\neg A), B_2{}^{g_1},\ldots, B_k{}^{g_1} \vdash (\neg(\neg B_1))$, I¬ from 1. and 3.
5. $((\neg A), B_2{}^{g_1},\ldots, B_k{}^{g_1} \vdash (\neg B_1)$, I¬ from 2. and 3.
6. $B_2{}^{g_1},\ldots, B_k{}^{g_1} \vdash (\neg(\neg A))$, I¬ from 4. and 5.
7. $B_2{}^{g_1},\ldots, B_k{}^{g_1} \vdash A$, R¬ from 6.

It is deduction of sequence

$$B_2{}^{g_1},\ldots, B_k{}^{g_1} \vdash A.$$

.

---

[§]Here $\vdash$ denotes $\rightarrow$.

This sequence is obtained from sequence (4) by deletion of first sentence from the hypothesizes list.

All rest hypothesizes are deleted from this list in the similar way.

Final sentence is the following:

$$\vdash A.$$

**Lm. 3:** Every natural propositional deduction of a sequence $\Gamma \rightarrow A$ satisfy the following condition: if $A$ is not true then there exists a sentence $C$ such that $C \in \Gamma$ and $C$ is not true.

**Proof of Lm. 3**: is realized by a recursion on number of sequences in the deduction of $\Gamma \rightarrow A$:

**1. Basis of recursion:** Let the deduction of $\Gamma \rightarrow A$ contains 1 sequence.

In that case a form of this sequence is $A \rightarrow A$ in accordance with the propositional natural deduction definition. Hence in this case the lemma holds true.

**2. Step of recursion:** *The recursion assumption:* Let the lemma holds true for every deduction containing no more than $n$ sequences.

Let the deduction of $\Gamma \rightarrow A$ contains $n + 1$ sequences.

In that case either this sequence is a NPC-axiom or $\Gamma \rightarrow A$ is obtained from previous sequences by one of deduction rules.

If $\Gamma \rightarrow A$ is a NPC-axiom then the proof is the same as for the recursion basis.

e) Let $\Gamma \rightarrow A$ be obtained from a previous sequence by $R$&.

In that case a form of this previous sequence is either the following $\Gamma \rightarrow (A\&B)$ or is the following $\Gamma \rightarrow (B\&A)$ in accordance with the definition of deduction. The deduction of this sequence contains no more than $n$ elements. Hence the lemma holds true for this deduction in accordance with the recursion assumption. If $A$ is not true then $(A\&B)$ is not true and $(B\&A)$ is not true. Hence there exists a sentence $C$ such that $C\in \Gamma$ and $C$ is not true in accordance with the lemma.

Hence in that case the lemma holds true for the deduction of sequence $\Gamma \rightarrow A$.

f) Let $\Gamma \rightarrow A$ be obtained from previous sequences by I&.

In that case forms of these previous sequences are $\Gamma_1 \rightarrow B$ and $\Gamma_2 \rightarrow G$ with $\Gamma = \Gamma_1,\Gamma_2$ and $A = (B\&G)$ in accordance with the definition of deduction. The lemma holds true for deductions of sequences $\Gamma_1 \rightarrow B$ and $\Gamma_2 \rightarrow G$ in accordance with the recursion assumption because these deductions contain no more than $n$ elements. In that case if $A$ is not true then $B$ is not true or $G$ is not true. Hence there exist a sentence $C$ such that $C$ is not true and $C \in \Gamma_1$ or $C \in \Gamma_2$.

Hence in that case the lemma holds true for the deduction of sequence $\Gamma \rightarrow A$.

g) Let $\Gamma \rightarrow A$ be obtained from a previous sequence by R¬.

In that case a form of this previous sequence is the following: $\Gamma \rightarrow (\neg(\neg A))$ in accordance with the definition of deduction. The lemma holds true for the deduction of this sequence in accordance with the recursion assumption because this deduction contains no more than $n$ elements. If $A$ is not true then $(\neg (\neg A))$ is not true. Hence there exists a sentence $C$ such that $C \in \Gamma$ and $C$ is not true.

Hence the lemma holds true for the deduction of sequence $\Gamma \to A$.

h) Let $\Gamma \to A$ be obtained from previous sequences by I¬.

In that case forms of these previous sequences are $\Gamma_1,G \to B$ and $\Gamma_2,G \to (\neg B)$ with $\Gamma = \Gamma_1, \Gamma_2$, and $A = (\neg G)$ in accordance with the definition of deduction. The lemma holds true for the deductions of sequences $\Gamma_1,G \to B$ and $\Gamma_2,G \to (\neg B)$ in accordance with the recursion assumption because these deductions contain no more than *n* elements.

If *A* is not true then *G* is true.

Either *B* is not true or $(\neg B)$ is not true. Hence there exists a sentence *C* such that either $C \in \Gamma_1,G$ or $C \in \Gamma_2,G$ and *C* is not true in accordance with the recursion assumption.

Hence in that case the lemma holds true for the deduction of sequence $\Gamma \to A$.

**The recursion step conclusion:** Therefore, in each possible case, if the lemma holds true for a deduction contained no more than *n* elements then the lemma holds true for a deduction contained $n + 1$ elements.

**The recursion conclusion:** Therefore the lemma holds true for a deduction of any length.

**Proof of Th. 3:** If a sentence *A* is naturally propositionally proven then there exists a natural propositional deduction of form $\to A$ (deduction from the empty list of hypothesizes). Hence, *A* is true in accordance with Lm.3.

**Proof of Th. 4:** Each tautology is naturally propositionally proven sentence by Th. 2. Each naturally propositionally proven sentence is a true sentence by Th.3. Therefore, every tautology is the true sentence.

**Proof of Th. 5:** If B∈[s]$(r,k)$ then ȸ(B) = $p^k(1-p)^{r-k}$ in accordance with Def. 22 and with (p10).

Since [*s*]$(r,k)$ contains $r!/(k!(r-k)!)$ elements then this theorem hold true according with (p9), (p10), and (p11).

**Proof of Th. 6:** In accordance with Def. 26: there exist natural numbers *n* and *k* such that $k-1 < a \le k$ and $k + n \le b < k + n + 1$, and Ŧ[$s_A$]$(r,a,b)$ := Ŧ[$s_A$]$(r,k,k+n)$.
The recursion on *n*:
**Basis of recursion:** Let $n = 0$.
In that case according Def. 25 and Def. 24:

Ŧ[$s_A$]$(r,k,k)$ = ŧ[$s_A$]$(r,k)$ = °«$\nu_r(s_A) = k/r$».

**Step of recursion:**
<u>*The recursion assumption:*</u> Let

Ŧ[$s_A$]$(r,k,k+n)$ = °«$k/r \le \nu_r[s_A] \le (k+n)/r$».

According to Def. 26:

Ŧ[$s_A$]$(r,k,k+n+1)$ = Ŧ[$s_A$]$(r,k,k+n)$ + ŧ[$s_A$]$(r,k+n+1)$.

According to the recursion assumption and according to Def. 25:

Ŧ[$s_A$]($r$,$k$,$k$ + $n$ + 1) = (°«$k/r \leq \nu_r[s_A] \leq (k+n)/r$» + °«$\nu_r(s_A) = (k+n+1)/r$»).

Hence according to Def.15:

Ŧ[$s_A$]($r$,$k$,$k$ + $n$ + 1) = (°«$k/r \leq \nu_r[s_A] \leq (k+n+1)/r$».

**The recursion step conclusion:** Therefore, if this theorem holds true for $n$ then one holds true for $n$ + 1.

**The recursion conclusion:** Therefore, this theorem holds true for any $n$.

**Proof of Th.7:** It follows from Th.5 and (p9) at once.

**Proof of Th. 8:** Because

$$\sum_{k=0}^{r}(k-rp)^2\cdot\frac{r!}{k!\,(r-k)!}p^k(1-p)^{r-k}=r\cdot p\cdot(1-p)$$

then if

$$J := \{k \in \mathbf{N} \mid 0 \leq k \leq r\cdot(p-\varepsilon)\} \cup \{k \in \mathbf{N} \mid r\cdot(p+\varepsilon) \leq k \leq r\}$$

then

$$\sum_{k\in J}\frac{r!}{k!\,(r-k)!}p^k(1-p)^{r-k} \leq \frac{p\cdot(1-p)}{r\cdot\varepsilon^2}.$$

Hence this theorem holds true according to (p3).

## *References*

1. Lyndon, R. *Notes on logic*, D. VAN NOSTRAND COMPANY, INC., 1966
2. Bernoulli, J. *Ars Conjectandi*, BASILEA, Impenfis THURNISORUM, Fratrum, 1713
3. Tarski, A. The Semantic Conception of Truth and the Foundations of Semantics, *Philosophy and Phenomenological Research*, **4**, 1944.
4. Mendelson, E. *Introduction to Mathematical Logic*, Chapmen&Hall, 1997, pp.11-50.
5. Gentzen, G., Untersuchungen Uber das Logische Schliessen, *Math. Z.* **39**, (1934), pp. 176-210, 405-431
6. Väth, M. *Nonstandard Analysis*, Birkhдuser Verlag, Basel, Boston, Berlin, 2007.
7. Quznetsov, G. Logical Foundation of Theoretical Physics, Nova Sci. Publ., NY, 2006.